\documentclass[a4paper, 12pt]{article}
\usepackage{enumerate, theorem}
\usepackage{amsmath, amsfonts, amssymb}
\usepackage[height=22.5cm, width=15cm]{geometry}
\usepackage[all]{xy}
\usepackage{mathrsfs, yfonts}

\DeclareMathOperator{\Sym}{Sym}

 \def\C{{\mathbb C}}

\newcommand{\parag}[1]{\paragraph{\sc{#1.}} }

\newtheorem{thm}{Theorem}[subsection]
\newtheorem{defn}[thm]{Definition}
\newtheorem{cor}[thm]{Corollaire}
\newtheorem{prop}[thm]{Proposition}
\newtheorem{lemma}[thm]{Lemma}

\begin{document}

\date{14/09/16}

\author{Daniel Barlet\footnote{Institut Elie Cartan : Alg\`{e}bre et G\'eom\`{e}trie,  \newline
Universit\'e de Lorraine, CNRS UMR 7502  and  Institut Universitaire de France.}.}

\title{Note on the quasi-proper direct image with value in a Banach analytic set.}

\maketitle

\parag{Abstract} We give a rather simple proof of the generalization of Kuhlmann's quasi-proper direct image theorem to the case of a map with values in  a Banach analytic set. The proof uses a generalization of the Remmert-Stein's theorem to this context.

\parag{AMS Classification} 32 H 02 - 32 K 05 - 32 C 25.  

\parag{Key words} Quasi-proper map. Quasi-proper direct image. Banach analytic set. 

\tableofcontents

\newpage

\section*{Introduction}
The aim of the present Note is to give  a rather simple proof for the generalization of  Kuhlmann's quasi-proper direct image theorem to the case of a map with values in Banach analytic set.  N. Kuhlmann  (and D. Mathieu in the Banach case, see [K.64], [K.66] and [M.00]) proved this direct image result for a ``semi-proper''  holomorphic map; this is a weaker hypothesis than quasi-proper, so their result is better  than the one presented here. But the quasi-proper case fits well with the situation we are mainly interested in the study of f-analytic family of cycles (see  [B.08], [B.13] and [B.15]\footnote{Nevertheless we prove the semi-proper case of the direct image theorem with values in $\mathcal{C}_{n}^{f}(M)$, the space of finite type $n-$dimensional cycles in the complex space $M$ in [B.15]  th. 2.3.2.}). Our simpler argument via the generalization of the Remmert-Stein's theorem does not work for a semi-proper map (see the remark following proposition \ref{Grosses fibres}). In the appendix we give an easy proof of the Remmert's direct image theorem in the proper finite case (with values in an open set of a Banach space) which is used in our proof  to make this Note  self-contained modulo the theorem  III 7.4.1. in [B-M 1].

\section{A simple  proof of Kuhlmann's quasi-proper \\ direct image theorem}

In order to show the strategy of proof for the generalization of Kuhlmann's theorem with values in an open set of a Banach space, we shall begin by a simple proof of the finite dimensional case using the  ``usual''  Remmert-Stein's theorem. 

\begin{thm}\label{Quasi-proper direct image}
Let \ $f : M \to N$ \ a quasi-proper holomorphic map between irreducible complex spaces. Then \ $f(M)$ \ is a closed analytic subset of \ $N$.
\end{thm}

\parag{Proof} The first point is to prove that \ $f(M)$ \ is closed in \ $N$. Let \ $y_{n} = f(x_{n})$ \ be a sequence in \ $f(M)$ \ and suppose that \ $y_{n}$ \ converges to a point \ $y \in N$. Let \ $V$ \ be a neighbourhood of \ $y$ \ in \ $N$ \ and \ $K$ \ a compact subset in \ $M$ \ such that for any \ $z \in V$ \ each irreducible component of \ $f^{-1}(z)$ \ meets \ $K$. Then we may assume that for \ $n$ \ large enough \ $y_{n}$ \ lies in \ $V$ \ and so that \ $x_{n}$ \ is chosen in \ $K$. So, up to pass to a sub-sequence we may assume that the sequence \ $(x_{n})$ \ converges to \ $x \in K$. Then the continuity of \ $f$ \ implies hat \ $y = f(x)$. So \ $f(M)$ \ is closed.\\
Consider now the integer \ $p : = \max\{ \dim_{x}f^{-1}(y), y \in N\}$ \ and define the set 
$$Z : = \{y \in N \ / \ \dim(f^{-1}(y)) = p \}.$$
 We want to show that \ $Z$ \ is a closed analytic subset in \ $N$. To show that \ $Z$ \ is closed consider a sequence \ $(y_{n})$ \ in \ $Z$ \ converging to a point \ $y$ \ in \ $N$. Let \ $V$ \ be a neighbourhood of \ $y$ \ in \ $N$ \ and \ $K$ \ a compact subset in \ $M$ \ such that for any \ $z \in V$ \ each irreducible component of \ $f^{-1}(z)$ \ meets \ $K$. Then for each \ $n$ \ large enough in order that \ $y_{n}$ \ lies in \ $V$, we may find a point \ $x_{n}\in K$ \ such that \ $f(x_{n}) = y_{n}$ \ and with \ $ \dim_{x_{n}}f^{-1}(y_{n}) = p$. So, up to pass to a subsequence we may assume that the sequence \ $x_{n}$ \ converges to \ $x \in K$. We have \ $f(x) = y$ ; consider now a \ $p-$scale \ $E :  = (U,B,j)$ \  around \ $x$ \ adapted to \ $f^{-1}(f(x))$. Then the map \ $g : = f \times(p_{U}\circ j) : j^{-1}(U\times B) \to W\times U$ \ is proper and finite if \ $W \subset V$ \ is an open neighbourhood of \ $x$ \ which is small enough. For each \ $x_{n} \in W$ \ the image of \ $g$ \ contains \ $\{x_{n}\}\times U$. So it contains also \ $\{x\}\times U$ \ and the fibre \ $f^{-1}(y)$ \ has dimension \ $p$ \ in \ $W$. So \ $y$ \ is in \ $Z$ \ and \ $Z$ \ is closed.\\
To prove the analyticity of \ $Z$ \ we argue as above and we remark that for any holomorphic function \ $h : U \times B \to \C^{m}$ \ with zero set \ $j(W)\cap (U \times B)$ \ we have \ $h(y,t) \equiv 0$ \ for \ $t \in U$ \ when \ $y $ \ is in \ $Z$. This will give (infinitely many)  holomorphic equations for \ $Z$ \ in \ $W$ \ using a Taylor expansion in \ $t$ \ of the holomorphic function \ $h$. This proves  the analyticity of $Z$.
Remark that \ $Z$ \ has dimension \ $ \leq m -p$ \ where $m := \dim M$.\\

Now we shall prove the theorem by induction on the integer \ $p \geq 0$. Note that, for \ $p = 0$ \ the result is already proved. So let assume that the result is proved for \ $p = q -1 $ \ and we shall prove it for \ $p= q$. If the closed analytic  set \ $Z$ \ introduced above is equal to \ $N$, the result is clear because \ $f(M) = N$. If not \ $N \setminus Z$ \ is an open dense set in \ $N$ \ because \ $N$ \ is irreducible and also \ $M \setminus f^{-1}(Z)$ \ is also an open dense set in \ $M$. Moreover the map  \ $\tilde{f} : M \setminus f^{-1}(Z) \to N \setminus Z$ \ induced by \ $f$ \ is again quasi-proper. The induction hypothesis gives that \ $f(M \setminus f^{-1}(Z))$ \ is a closed analytic subset in \ $N \setminus Z$ \ and it has dimension bigger or equal to \ $m-(p-1) = m - p + 1$. As \ $Z$ \ has dimension \ $ \leq m - p$ \ we may apply Remmert-Stein theorem to obtain that the closure of \ $f(M \setminus f^{-1}(Z))$ \ in \ $N$ \ is again an analytic set. This conclude the proof as we have \ $f(M) = \overline{f(M \setminus f^{-1}(Z))} \cup Z$, because we know that \ $f(M)$ \ is closed.$ \hfill \blacksquare$

\parag{Remarks}
\begin{enumerate}
\item The semi-properness of \ $f$ \ would be enough to get the fact that \ $f(M)$ \ is closed.
\item The restriction of a quasi-proper map to a  closed saturated analytic subset $X$ in $M$\footnote{A subset $X$ of $M$  is saturated for $f$  if $x \in X$ implies $f^{-1}(x) \subset X$.} is again quasi-proper. This is also true for a semi-proper map. 
\item When $f$ is quasi-proper to have quasi-properness for $f$ restricted to closed analytic subset $X \subset M$ it is enough for $X$ to be an  union of some  irreducible components of  $f^{-1}(y)$ for  $y \in f(X)$. 
\item This property of $X$  is not enough, in general, in the case of a semi-proper map $f$. This is precisely what happen for the subset $Z$ introduced below or in the proposition \ref{Grosses fibres}.
\end{enumerate}

\newpage

\section{The case with values in a Banach analytic set}

\begin{defn}\label{q-proper}
Let \ $f : M \to S$ \ be a holomorphic map from a  reduced complex space \ $M$ \ to a Banach analytic set \ $S$. We say that {\bf \ $f$ \ is quasi-proper} \ when for each point \ $s \in S$ \ there exists a neighbourhood \ $V$ \ of \ $s$ \ in \ $S$ \ and a compact set \ $K$ \ in \ $M$ \ such that,  for {\bf each} \ $\sigma \in V$, {\bf each} irreducible component of \ $f^{-1}(\sigma)$ \ meets \ $K$.
\end{defn}

Note that \ $f$ \ quasi-proper implies semi-proper\footnote{Recall that a continuous map $f : M \to S$ between Hausdorff topological spaces is semi-proper if for any $s \in S$ there exists a neighbourhood $V$ of $s$ in $S$ and a compact set $K$ in $M$ such that for any $\sigma \in V$ the fibre $f^{-1}(\sigma)$ meets $K$. This is a topological notion contrary to the quasi-properness which asks that $M$ is a complex space and that the fibres of $f$ are (closed) analytic subsets in $M$ (but no complex structure in needed on $S$).} as the definition above implies that \ $f(K) \cap V = f(M)\cap V$. In particular this condition implies that \ $f(M)$ \ is closed in \ $S$. So our result is local on $S$ and  it is enough to prove the Kulmann's theorem with value in an open set of a Banach space.

\begin{prop}\label{Grosses fibres}
Let \ $f : M \to S$ \ be a holomorphic map from a  reduced and irreducible complex space \ $M$ \ to a Banach analytic set \ $S$. We assume that \ $f$ \ is quasi-proper. Let \ $p : = \max\{\dim f^{-1}(s), s \in S\}$, and define 
 $$ Z : = \{s \in S\ / \ \dim f^{-1}(s) = p \}.$$
  Then \ $Z$ \ is a closed analytic subset in \ $S$ \ which is finite dimensional (so \ $Z$ \ is a reduced complex space embedded in \ $S$; see the theorem III 7.4.1 de [B-M 1]).
\end{prop}

\parag{remark} This result is not true in general for a semi-proper  holomorphic map:\\
Let $\pi : M \to N$ be an infinite connected cover of a complex manifold $N$ of dimension $n \geq 2$. Let $y \in N$ and let $V$ be a relatively compact open set containing $y$ such that $\pi$ admit a continuous section $\sigma : \bar V \to \bar W$ where $ W$ is a relatively compact open set in $M$. Let $(y_{\nu})$ be a sequence of points in $V\setminus \{y\}$ with limit $y$ and let $(x_{\nu})$ be a discrete sequence in $M \setminus W$ such that $\pi(x_{\nu}) = y_{\nu}$ for each $\nu \in \mathbb{N}$. Let $ \tau : \tilde{M} \to M$ be the blow-up of $M$ at each point $x_{\nu}$ and put $\tilde{\pi} := \tau\circ \pi$. Then $\tilde{\pi}$ is semi-proper and the subset $Z$ of $M$ where the fibre of $\tilde{\pi}$ has dimension $n-1$ is exactely the subset $\{y_{\nu}, \nu \in \mathbb{N}\}$ which is not closed in $N$.\\

\parag{First step of the proof:  $Z$ \ is closed in \ $S$} Let \ $(s_{\nu})_{\nu \geq 0}$ \ be a sequence in \ $Z$ \ which converges to a point \ $s $ \ in \ $S$. Let \ $V$ \ be an open neighbourhood of \ $s$ \ and \ $K$ \ a compact subset of \ $M$ \ such that for any \ $\sigma \in V$ \ each  irreducible component of \ $f^{-1}(\sigma)$ \ meets \ $K$. For \ $\nu \gg 1$ \ we have \ $s_{\nu} \in V$, and if \ $\Gamma_{\nu}$ \ is a \ $p-$dimensional irreducible component of \ $f^{-1}(s_{\nu})$, the intersection \ $\Gamma_{\nu}\cap K$ \ is not empty, and we may choose some \ $x_{\nu} \in \Gamma_{\nu}\cap K$. Up to pass to a sub-sequence, we may also assume that the sequence \ $(x_{\nu})$ \ converges to a point \ $x \in K$. Of course we shall have \ $x \in f^{-1}(s)$. Choose now a \ $p-$scale \ $E : = (U,B,j)$ \ such that \ $x $ \ is in \ $j^{-1}(U\times B)$ \ and such that \ $j^{-1}(\bar U\times \partial B) \cap f^{-1}(s) = \emptyset$. This is possible because we know that \ $\dim f^{-1}(s) \leq p$. Then, up to shrink \ $V$ \ around \ $s$, we may assume that for any \ $\sigma \in V$ \ we have \ $f^{-1}(\sigma) \cap j^{-1}(\bar U\times \partial B)= \emptyset$. This means that for \ $\nu \gg 1$ \ the scale \ $E$ \ is adapted to \ $\Gamma_{\nu}$ \ and with \ $\deg_{E}(\Gamma_{\nu}) = k_{\nu} \geq 1$.\\
 Define on the open set $j^{-1}(U\times B)\cap f^{-1}(V) \subset M$ the holomorphc map
  $$g := (p_{U}\circ j)\times f : j^{-1}(U\times B)\cap f^{-1}(V) \to U\times V $$
   where $p_{U} : U \times B \to U$ is the projection and where $V$ is the open set in $S$ defined by the condition $s \not\in f(j^{-1}(\bar U\times \partial B))$. 
\parag{Second step: $g$ is a closed map with finite fibres} The finiteness of fibres is obvious because for \ $\sigma \in V$ the intersection \ $f^{-1}(\sigma)\cap j^{-1}( \bar U\times \partial B) $ is empty, so \ $ j(g^{-1}(t, \sigma)) \subset \{t\}\times B$ \ is a compact analytic subset in a polydisc, so a finite set. \\
To show the closeness of $g$, choose a closed set $F$ in $f^{-1}(V)\cap j^{-1}(U\times B)$ and   a sequence $(t_{\nu},\sigma_{\nu})$ \ in \ $g(F) \cap (U\times V)$ \ converging to a point \ $(t, \sigma) \in U \times V$, and let \ $(x_{\nu})$ \ be a sequence in \ $ F$ \ such that \ $g(x_{\nu}) = (t_{\nu}, \sigma_{\nu})$. As \ $j^{-1}(\bar U \times \bar B)$ \ is compact, we may assume, up to pass to a sub-sequence, that the sequence \ $(x_{\nu})$ \ converges to some \ $x \in j^{-1}(\bar U \times \bar B) \subset M$. But the limit of  \ $g(x_{\nu}) = (t_{\nu}, \sigma_{\nu})$ \ is in $U\times V$ by assumption and this implies that  $x$ is in \ $j^{-1}(U\times \bar B)$. As it cannot be in \ $j^{-1}(U \times \partial B)$ because $f(x)$ is in $V$, $x$ \ lies in \ $f^{-1}(V)\cap j^{-1}(U\times B)$  and so $x$ is in $F$. This proves the closeness of \ $g$. \\
Now the Remmert's theorem in the proper finite case, but with values in a Banach analytic set\footnote{So we have to use the proper case with finite fibres of the Remmert' theorem with values in an open set of a complex Banach space to prove this proposition. To make this Note self-contained (modulo the theorem  III 7.4.1. in [B-M 1]) we give a simple proof of this case in the appendix (section 4). See also [B-M 1] chapter III for the general case.}, applies and shows that the image of \ $g$ \ is a reduced complex space. In this case it is clear that the cardinal of the fibres is locally bounded, so \ $k_{\nu}$ \ is bounded up to shrink \ $U$ \ and \ $V$.\\
So, up to pass to a subsequence, we may assume that \ $k_{\nu}$ \ is constant equal to \ $k \geq 1$. Then it is easy to see, again up to pass to a sub-sequence and to shrink \ $U$, that the sequence \ $j(\Gamma_{\nu})$ \ converges to a multiform graph. This implies that \ $f^{-1}(s)$ \ has dimension \ $p$, so \ $s$ \ is in \ $Z$.\\
We shall denote $\Xi \subset (U \times V)$ the image of $g$.

\parag{Third step}  We shall prove now that, assuming that we choose \ $V$ \ small enough around the given point $s$ in $Z$, the set
$$ Z' : = \{ \sigma \in V \ / \  \{\sigma\}\times U \subset \Xi \} $$
is a closed analytic subset in \ $V$. Remark that $Z'$ is contained in $Z$ but it may be smaller that $Z$ because the fibre of a point $s \in V \setminus Z'$ may be of dimension $p$ via a component of $f^{-1}(s)$ which does not meet $j^{-1}(U\times B)$.\\
 Up to shrink \ $V$ \ and \ $U$ \ we may assume that \ $V$ \ is a closed  analytic subset of an open set \ $\mathcal{U}$ \ in some Banach space, and  that we have a holomorphic map
$$ \Phi : \mathcal{U} \to Hol(\bar U, F) $$
where \ $F$ \ is a Banach space, such that the associated holomorphic map  $$ \tilde{\Phi} : \mathcal{U}\times U \to F$$
 satisfies \ $  \tilde{\Phi}^{-1}(0) \cap (V \times U) = \Xi$. Then it is clear that we have \ $Z' = \Phi^{-1}(0)$.\\
Now remark that \ $Z'\times U$ \ is a closed analytic subset of \ $\Xi$ \ which is finite dimensional. So this implies that \ $Z'$ \ is also finite dimensional and we have \ $\dim Z' \leq \dim M - p$.\\
If we cover the compact set \ $f^{-1}(s) \cap K$ \ by finitely many $p-$scales as above, we obtain that \ $Z$ \ is locally a finite union of such $Z'$ as above and then $Z$  is a finite dimensional analytic set of dimension \ $\leq \dim M - p$ \ near the point \ $s$;  and, as we already know that \ $Z$ \ is closed, it is  a closed finite dimensional analytic set of  dimension \ $\leq \dim M - p$ \ in $S$.$\hfill \blacksquare$

\begin{thm}\label{Kuhl-banach}
Let  $f : M \to N$  a quasi-proper holomorphic map between an irreducible complex space $M$ and a Banach analytic set $S$. Then f(M) is a closed analytic subset in S which is locally finite dimensional .
\end{thm}

\parag{Proof} We shall prove the theorem by induction on the integer \ $p$ defined in the previous proposition. Note that for a given \ $M$, we have always \ $p \leq \dim M < + \infty$ \ as we assume \ $M$ irreducible.\\
The case \ $p = 0$ \ is clear because in this case we have \ $Z = f(M)$.\\
Assume the theorem proved when \ $ p \leq q-1$ \ for some \ $q \geq 1$ \ and we shall prove it for \ $p = q$. \\
From the previous steps we know that \ $Z$ \ is a closed analytic subset in \ $S$ \ and that it has finite dimension \ $\leq \dim M - p$. Let us consider now the map 
$$ \tilde{f} : M \setminus f^{-1}(Z) \to  S \setminus Z $$
induced by \ $f$.
It is clearly quasi-proper and for this map we have \ $\tilde{p} \leq q - 1$. So, by the induction hypothesis, the image of \ $\tilde{f}$ \ is a closed analytic subset in \ $S\setminus Z$ \ which is irreducible of finite dimension 
  $$\dim Im(\tilde{f}) \geq \dim M -  \tilde{p}  \geq \dim M - q + 1 \geq \dim Z + 1.$$
Now we want to apply the Remmert-Stein theorem to conclude. This is clear in the case where \ $S$ \ is a finite dimensional complex space because the dimensions satisfy the desired inequality, but we want now to treat the case where \ $S$ \ is a Banach analytic set. As the problem is local, we may replace \ $S$ \ by a open set in a Banach space, and we apply the  generalization of Remmert-Stein theorem obtained in the next section to conclude.$\hfill \blacksquare$

\newpage

\section{The Remmert-Stein' theorem in a Banach space}

\begin{thm}\label{R-S Banach}
Let \ $\mathcal{U}$ \ an open set in a Banach space and let \ $A \subset \mathcal{U}$ \ be a closed analytic susbset of dimension\ $d$. Let \ $X \subset \mathcal{U}\setminus A$ \ a closed irreducible analytic subset in \ $\mathcal{U}\setminus Z$ \ of finite dimension \ $\geq d +1$. Then \ $\bar X$ \ is a closed irreducible analytic subset in \ $\mathcal{U}$ \ of finite dimension equal to \ $\dim X$.
\end{thm}

Here is a first step of the proof of this theorem.

\begin{lemma}\label{cut-discret}
Let \ $X \subset \mathcal{U}$ \ be a closed analytic subset in an open set \ $\mathcal{U}$ \ in a Banach space \ $E$. Assume that \ $X$ \ has finite  pure dimension \ $\leq n$ \ and that \ $X$ \ is  countable at infinity. Let \ $a \in \mathcal{U} \setminus X$ \ and \ $x \in X$. Assume also that a \ $n-$dimensional linear subspace \ $L$ \ containing \ $a$ \ is given. Then there exists a linear closed codimension \ $n$ \ subspace \ $P$ \ in \ $E$ \ containing \ $a$ \ and \ $x$,  such \ $P$ \ is transversal to \ $L$ \ at \ $a$ \ and  that the set \ $X \cap P$ \ is discrete and countable.
\end{lemma}

\parag{Proof} We shall make an induction on \ $n \geq 0$. The case \ $n = 0$ \ is trivial. So let assume that the case \ $n-1$ \ is proved and consider \ $X$ \ of pure dimension \ $d \leq n$. Choose for each irreducible component \ $X_{\nu}$ \ of \ $X$ \ a point \ $x_{\nu} \in X_{\nu}$ \ which is not equal to \ $x$\footnote{As \ $n \geq 1$ \ an irreducible component is not equal to \ $\{x\}$.}. Now choose an hyperplane \ $H$ \ containing \ $a$ \ and \ $x$ \ and such that \ $x_{\nu}$ \ is not in \ $H$ \ for each \ $\nu$, and that \ $H$ \ is transversal to \ $L$ \ at \ $a$. Such a hyperplane exists thanks to Baire's theorem. Then \ $X\cap H$ \ is purely \ $(d-1)$-dimensional and contains \ $x$. The induction hypothesis gives us a co-dimension \ $n-1$ \ subspace \ $\Pi$ \ containing \ $a$ \ and \ $x$, transversal at \ $a$ \ to \ $L \cap H$,  such that \ $\Pi\cap (X\cap H)$ \ is discrete and countable. Now  the co-dimension of \ $ P : = \Pi \cap H$ \ is exactly equal to \ $n$ and is transversal to $L$ at $a$. This completes the proof. $\hfill \blacksquare$

\parag{Proof of the theorem \ref{R-S Banach}} We follow the proof of  the proposition II 4.8.3 case (i) given in [B-M 1]  but with an ambiant  infinite dimensional complex Banach space. \\
The theorem \ref{R-S Banach} follows from this result as in the finite dimensional case.\\

First, the following point is not completely obvious : the local compactness of \ $X \cup A$. \\
Of course \ $A$ \ and \ $X$ \ are locally compact by assumption. So it is enough to show that this union is locally compact near a point \ $a \in A$. Remark first that we have
$$ \partial X \subset \partial \mathcal{U} \cup A.$$
and $X \cup A$ is closed in $\mathcal{U}$. So if we consider \ $\varepsilon > 0 $ \ small enough, we have \ $\bar B(a,2\varepsilon) \subset \mathcal{U}$ \ and \ $\bar B(a, \varepsilon) \cap A$ \ is a compact set. Now we want to show that \ $\bar B(a,\varepsilon) \cap (X \cup A)$ \ is compact. So it is enough to show that any sequence \ $(x_{\nu}) $ \ in \ $\bar B(a,\varepsilon) \cap X$ \ admits a converging sub-sequence to a point in $\bar B(a, \varepsilon)$. Let
$$\delta_{\nu} : = d(x_{\nu}, A \cap \bar B(a,\varepsilon)).$$
 If \ $ \delta_{\nu}$ goes to  $0$ \ for some sub-sequence, then there is a  Cauchy sub-sequence converging to some \ $a' \in \bar B(a,\varepsilon)\cap (X \cup A)$. If \ $\delta_{\nu}\geq  \alpha > 0$ \ for all large enough \ $\nu$, then the sequence is contained in the subset
$$\Big(\bar B(a,\varepsilon) \setminus \{ x \ / \ d(x, A\cap\bar B(a,\varepsilon))\geq \alpha)\} \Big) \cap X.$$
This subset of \ $X$ \ is compact, because it is closed and a sequence in it cannot approach neither \ $\partial \mathcal{U}$ \ (as $\bar B(a, 2\varepsilon)\subset \mathcal{U}$, any point $z \in \bar B(a, \varepsilon)$ satisfies $d(z, \partial\mathcal{U}) \geq \varepsilon$) nor \ $A$, and \ $X$ \ is locally compact (so any discrete sequence has to approach the boundary of $X$).\\

Assume that $X$ has  dimension $d+k$ with $k \geq 1$. For \ $a \in A$ \ a fixed smooth point, let \ $V$ \ be a \ $d+k$ \ sub-manifold through \ $a$ \ containing \ $A$ (we use here the theorem III 7.4.1 in [B-M 1]) and  choose a co-dimension \ $d+k$ \ plane \ $P$ \ transversal to \ $V$ \ at \ $a$ \ meeting some point in \ $X \cap B(a,\varepsilon)$ \ and such that \ $P \cap X$ \ is discrete. \\
Up to shrink $V$ and taking a small ball $\mathcal{B}_{0}$ with center $a$ in $P$, we may assume that an open neighbourhood of $a$ in $B(a, \varepsilon)$ is isomorphic to the product $V \times \mathcal{B}_{0}$. And from our construction we have $X \cap (\{a\} \times \mathcal{B}_{0})$ which is discrete.
\parag{Claim}  We can choose two small balls \ $\mathcal{B}' \subset \mathcal{B}$ \ in $P$ with center \ $a$ \ contained in \ $\mathcal{B}_{0} \subset B(a,\varepsilon)\cap P$ \ such that \ $(\mathcal{B}\setminus \mathcal{B}') \cap X = \emptyset $ \   and then we may find an open neighbourhood \ $U$ \ of \ $a$ \ in \ $V$ \ in order that the projection on $U$ of the set  \ $\bar X \cap (U \times \mathcal{B}') = (X \cup A) \cap (U \times \mathcal{B})$ \ is proper (this means closed with compact fibres) and  that its restriction to $X \cap ((U\setminus A)\times \mathcal{B})$ has finite  fibres.
\parag{proof} The choice of two arbitrary small balls \ $\mathcal{B}' \subset \mathcal{B}$ with the  first condition $(\mathcal{B}\setminus \mathcal{B}') \cap X = \emptyset $\  is easy because we know that
 $X \cap (\{a\} \times \mathcal{B}_{0})$ is discrete : the distances to $a$ of the points in $X \cap(\{a\}\times \mathcal{B}_{0})$ are in a discrete subset in $]0, \varepsilon[$ and we can choose $r' < r$ with $r$ arbitrary small such that $]r', r[$ avoids these values. But the subset $(X \cup A) \cap \bar B(a, \varepsilon)$ is compact and then, as $A \cap (\bar{\mathcal{B}} \setminus \mathcal{B}')$ is empty, the subset $\bar X \cap (\{a\}\times \mathcal{B})$ is compact. Then there exists an open neighbourhood $U$ of $a$ in $V$ such that the projection of $\bar X \cap (U \times \mathcal{B})$ on $U$ is proper :\\
Let $K$ be a compact neighbourhood in $\bar X \cap (V \times \mathcal{B})$ of the compact set $\bar X \cap (\{a\}\times \mathcal{B})$ (remember that we proved that $\bar X = X \cup A$ is locally compact). Then for $K$ small enough the distance of a point in $K$ to the closed set $\{a\}\times (\bar{\mathcal{B}}\setminus \mathcal{B}')$ is bounded below by a positive number, and so there exists an open neighbourhood $U$ of $a$ in $V$ such that $U \times (\bar{\mathcal{B}}\setminus \mathcal{B}')$ does not meet $K$. This is enough to prove the claim thanks to the following simple remark, as $\bar X \cap ((U\setminus A)\times \mathcal{B}) = X \cap ((U\setminus A)\times \mathcal{B})$ because $A \simeq A \times \{a\} \subset U \times \mathcal{B}'$ .

\parag{Remark} Any compact analytic subset in $X \cap B(a, \varepsilon)$ is finite: assume that $Z$ is a connected component of such a compact analytic subset. As compactness implies that there are at most finitely many $Z$, it is enough to prove that $Z$ has at most one point.\\
Assume that $z_{1}\not= z_{2}$ are two points in $Z$. Choose a continuous linear form on $E$ such that $l(z_{1}) \not= l(z_{2})$. As the map $l_{\vert Z}: Z \to \C$ is holomorphic, it has to be constant (maximum principle and connectness). Conclusion : $Z$ has at most one point.\\

The end of the proof of the theorem \ref{R-S Banach} is now analogous to the end of the proof in the finite dimensional setting (we are in the ``easy'' case where the dimension of $X$ is strictly bigger than $d$ the dimension of $A$ ; see the remark before the lemma 4.8.7 in  [B-M 1] ch.III).$\hfill \blacksquare$

\parag{Remarks}\begin{enumerate}
\item Along the same line, it is not difficult to prove the analog of the Remmert-Stein's theorem with equality of dimensions, assuming there exists an open set in $\mathcal{U}$  meeting each irreducible component of $A$ and in which $X \cup A$ is an analytic subset.
\item  In our proof we use in a crucial way the metric in the ambient Banach space. If the extension of this proof to an open set in a Frechet space may be easy, the case of a non metrizable e.l.c.s.s. (for instance a dual of Frechet space ) does not seem clear. Note that the proper case is proved with values in an open set of  any e.l.c.s.s.  in [B-M 1] ch.III.
\end{enumerate}

\section{Appendix: the  direct image theorem in the proper finite case}

\parag{Notation}An e.l.c.s.s. is a locally convex topological vector space which is separated and sequentially complete. For the definition we shall use of a holomorphic map on a reduced complex space $M$  with values in a e.l.c.s.s. see [B-M 1] ch. III  definition 7.1.1.\\
Note that the general proper case of this result is given in {\it loc. cit. th. III 7.4.3}\  but the proof of the proper finite case is much more simple.

To have a simple and self-contained  proof (modulo the theorem  III 7.4.1. in [B-M 1]) of the quasi-proper direct image theorem discussed in this Note, we shall give a easy proof of the direct image theorem in the proper finite case with values in an e.l.c.s.s.(this result has been used in the Banach case in our proof).\\
This allows to use the  proof of the previous section  also in the proper case, thanks to the following lemma.

\begin{lemma}\label{quasi-propre}
Let $f : M \to S$ a holomorphic map from a reduced complex space to an open set of an e.l.c.s.s. Assume that $f$ is closed with compact fibres (so $f$ is proper). Then $f$ is quasi-proper.
\end{lemma}

\parag{proof} If $f$ is not quasi-proper at $s \in S$ there exists, for each open neighbourhood $V_{\nu}$ of $s$ and for each compact $K$ in $M$ a point $\sigma_{\nu, K}\in V_{\nu}$ and an irreducible component $\Gamma_{\nu,K}$ of $f^{-1}(\sigma_{\nu, K})$ which does not meet $K$. Now, for a fixed $K$ we may choose the points $\sigma_{\nu, K}$ to be distinct. Let $W$ be a relatively compact open set in $M$ containing the compact fibre $f^{-1}(s)$ and let $K := \bar W$. Then for each $\nu$ we can find a point $x_{\nu}$ which is in $\Gamma_{\nu, K}$. So the point $x_{\nu}$ is in the closed set $M \setminus W$. Consider the subset $F := \{ x_{\nu} \}$ in $M \setminus W$. It is a closed set in $M\setminus W$ (and also in $M$) because if $x \in \bar F \setminus F$, the point $x$ is limit of an ultra-filter of points in $F$, and so $f(x)$ is equal to $s$, because the intersection of the $V_{\nu}$ for any ultra-filter is reduced to $\{s\}$. But as $x \in M \setminus W$ and $f^{-1}(s) \subset W$, this is a contradiction. So the map $f$ is quasi-proper at any point $s \in S$.$\hfill \blacksquare$

\begin{prop}\label{proper finite case}
Let $f : M \to S$ be a holomorphic map of a complex reduced space $M$ to an open set $S$ in a e.l.c.s.s. Assume that $f$ is closed with finite fibres (so proper and finite). Then $f(M)$ is a closed analytic subset which is locally of finite dimension in $S$.
\end{prop}

\parag{proof} As $f(M)$ is the locally finite union of the sets $f(M_{i}), i \in I$, where $M_{i}, i \in I$ is the set of irreducible components of $M$ because $f$ is quasi-proper and finite thanks to the previous lemma, we may assume that $M$ is irreducible. Let $m := \dim M$. As $f(M)$ is a closed set in $S$, it is enough to show that $f(M)$ is an finite dimensional analytic subset near a point $s = f(x)$ in $f(M)$.\\
 Let $G$ be the ambiant vector space of $S$. Choose a finite co-dimensional closed affine subspace $H$ containing $s$  in the e.l.c.s.s. $G$ such that the point $s$ is isolated in  $H \cap f(M)$ which is maximal for this property.\\
 Such a $H$ exists because it is easy to construct a sequence $(H_{q})$ of closed co-dimension $q$ affine subspaces in $G$ containing $s$ such that the germ of $f^{-1}(H_{q})$ at $f^{-1}(s)$ gives a strictly decreasing sequence of analytic germs as long as $f^{-1}(H_{q})$ is strictly bigger than the finite set $f^{-1}(s)$. Then at the last step $q_{0}$ we have the equality of germs $f^{-1}(H_{q_{0}}) = f^{-1}(s)$ in $M$.
 If $H_{q_{0}}$ is not maximal with this property, replace $H_{q_{0}}$ by a maximal closed affine subspace containing $s$ and containing $H_{q_{0}}$ with this property.\\
 Write $G = \C^{n}\oplus H^{0}$ where $H^{0}$ is the closed co-dimension vector subspace in $G$ directing $H$ and let $p : G \to \C^{n}$ and $p_{H} : G \to H$  be the projections.\\
 Then the germ of analytic map $g : (M, f^{-1}(s)) \to (\C^{n}, p(s))$ is finite by construction. So its image is an analytic germ. If it is not equal to $(\C^{n}, p(s))$ we can find a line $\Delta$ through $p(s)$ in $\C^{n}$ such that $g^{-1}(p(s)) = f^{-1}(s)$. Then the affine space  $p^{-1}(\Delta)$ contradicts the maximality of $H$. So the germ $g$ is surjective and we can find an open polydisc $U$ with center $p(s)$ in $\C^{n}$ such that $g : M' := g^{-1}(U) \to U$ is a $k-$sheeted branch covering. So we have a holomorphic map $ \varphi : U \to  \Sym^{k}(M')$ classifying the fibres of $g$. Composed with the holomorphic map induced by $\Sym^{k}(p_{H}\circ f)$ we obtain that $f(M')$ is a multiform graph of $U$ contained in $U \times H$ (see [B-M 1] chapter III section 7.2). As we have $f(M') = p^{-1}(U)\cap M$ in a neighbourhood of $s$, we conclude that $f(M)$ is a finite dimensional analytic subset in $G$ in a neighbourhood of $s$.$\hfill \blacksquare$

\newpage

\section{Bibliography}

\begin{itemize}

\item{[B.08]} Barlet, D. \textit{Reparam\'etrisation universelle de familles f-analytiques de cycles ...} Comment. Helv. 83 (2008), pp. 869-888.

\item{[B.13]} Barlet Daniel, {\it Quasi-proper meromorphic equivalence relations}, Math. Z. (2013), vol. 273, pp. 461-484.

\item{[B.15]} Barlet, D. {\it Strongly quasi-proper maps and the f-flattening theorem}, math-arXiv:1504.01579 (41 pages)

\item{[B-M 1]} Barlet, D. et Magn\'usson, J. {\it Cycles analytiques complexes. I. Th\'eor\`{e}mes de pr\'eparation des cycles}, Cours Sp\'ecialis\'es, 22. Soci\'et\'e Math\'ematique de France, Paris  (2014).

\item{[K.64]} Kuhlmann, N. {\it \"Uber holomorphe Abbildungen komplexer Ra\"ume} Arch. der Math. 15, (1964) pp. 81-90.

\item{[K.66]}  Kuhlmann, N. {\it Bemerkungen \"uber holomorphe Abbildungen komplexer Ra\"ume} Wiss. Abh. Arbeitsgemeinschaft Nordrhein-Westfalen 33, Festschr. Ged\"aachtnisfeier K. Weierstrass, (1966) pp.495-522.

\item{[M.00]} Mathieu, D. {\it  Universal reparametrisation of a family of cycles : a new approach to meromorphic equivalence relations} Ann. Inst. Fourier (Grenoble) vol. 50 fasc. 4, (2000) pp. 1155-1189.

\end{itemize}

\end{document}